\documentclass[leqno]{article}

\title{Simultaneous approximation by Bernstein polynomials with integer coefficients}
\author{Borislav R. Draganov}
\date{}

\usepackage{amsmath,amsthm,amsfonts}
\usepackage[mathscr]{eucal}
\usepackage[english]{babel}
\usepackage{fancyhdr}

\allowdisplaybreaks[1]

\newtheorem{thm}{Theorem}[section]

\newtheorem{lem}[thm]{Lemma}
\theoremstyle{remark}
\newtheorem{rem}[thm]{\bf Remark}

\numberwithin{equation}{section}

\newcommand{\thmref}[1]{Theorem~\ref{#1}}

\newcommand{\lemref}[1]{Lemma \ref{#1}}

\newcommand{\N}{\mathbb{N}}
\newcommand{\R}{\mathbb{R}}
\newcommand{\Z}{\mathbb{Z}}
\newcommand{\la}{\left\langle}
\newcommand{\ra}{\right\rangle}

\begin{document}

\maketitle
\bigskip

\thispagestyle{fancy}
\fancyhf{}
\renewcommand{\headrulewidth}{0pt}
\lhead{}
\renewcommand{\footrulewidth}{0.3pt}
\lfoot{\footnotesize This work was supported by grant DN 02/14 of the Fund for Scientific Research of the Bulgarian Ministry of Education and Science.}

\begin{abstract}
We prove that several forms of the Bernstein polynomials with integer coefficients possess the property of simultaneous approximation, that is, they approximate not only the function but also its derivatives. We establish direct estimates of the error of that approximation in uniform norm by means of moduli of smoothness. Moreover, we show that the sufficient conditions under which those estimates hold are also necessary.
\end{abstract}

\bigskip
\noindent
{\footnotesize \leftskip25pt \rightskip25pt {\sl  AMS} {\it classification}: 41A10, 41A25, 41A28, 41A29, 41A35, 41A36.\\[2pt]
\par}
\bigskip

\section{Main results}

The Bernstein operator or polynomial is defined for $f\in C[0,1]$ and $x\in [0,1]$ by
\[
B_n f(x):=\sum_{k=0}^n f\left(\frac{k}{n}\right)p_{n,k}(x),\quad p_{n,k}(x):=\binom{n}{k}x^k (1-x)^{n-k}.
\]
It is known that if $f\in C[0,1]$, then
\[
\lim_{n\to\infty} \|B_n f-f\|=0,
\]
where $\|\circ\|$ is the sup-norm on the interval $[0,1]$. A best possible estimate of that convergence can be given by the Ditzian-Totik modulus of smoothness $\omega_\varphi^2(f,t)$ of the second order with a varying step, controlled by the weight $\varphi(x):=\sqrt{x(1-x)}$, in the uniform norm on the interval $[0,1]$ (see \cite[Chapter 2]{Di-To:Mod}). For all $f\in C[0,1]$ and $n\in\N_+$ there holds (see \cite[Chapter 10, (7.3)]{De-Lo:CA}, or \cite[Theorem 6.1]{Bu:BP})
\begin{equation}\label{eqB}
\|B_n f-f\|\le c\,\omega_\varphi^2(f,n^{-1/2}).
\end{equation}
Above and henceforward $c$ denotes a positive constant, not necessarily the same at each occurrence, whose value is independent of $f$ and 
$n$. Instead of $\omega_\varphi^2(f,t)$ we can use the moduli defined in \cite{Iv:Dir,Iv:Char}, \cite{Dr-Iv:Char}, or 
\cite{Dz-Sh,Ko-Le-Sh:2,Ko-Le-Sh,Ko-Le-Sh:3}.

Kantorovich \cite{Ka} (or e.g.\ \cite[pp.\ 3--4]{Fe}, or \cite[Chapter 2, Theorem 4.1]{Lo-Go-Ma:CA}) introduced an integer modification of 
$B_n$. It is given by
\[
\widetilde{B}_n(f)(x):=\sum_{k=0}^n \left[f\left(\frac{k}{n}\right) \binom{n}{k}\right] x^k (1-x)^{n-k}.
\]
Above $[\alpha]$ denotes the largest integer that is less than or equal to the real $\alpha$. L.~Kantorovich showed that if $f\in C[0,1]$ is such that 
$f(0),f(1)\in\Z$, then
\begin{equation*}
\lim_{n\to\infty} \|\widetilde{B}_n (f)-f\|=0.
\end{equation*}
Clearly, the conditions $f(0),f(1)\in\Z$ are also necessary in order to have $\lim_{n\to\infty} \widetilde{B}_n (f)(0)=f(0)$ and 
$\lim_{n\to\infty} \widetilde{B}_n (f)(1)=f(1)$, respectively. 

Following L.~Kantorovich and applying \eqref{eqB}, we get a direct estimate of the error of $\widetilde B_n$ for $f\in C[0,1]$ with 
$f(0),f(1)\in\Z$. For $x\in [0,1]$ and $n\in\N$ we have
\begin{equation}\label{eq[B]}
\begin{split}
|\widetilde{B}_n (f)(x)-f(x)| &\le |B_n f(x)-f(x)| \\
& + \sum_{k=1}^{n-1} \left(f\left(\frac{k}{n}\right) \binom{n}{k}-\left[f\left(\frac{k}{n}\right) \binom{n}{k}\right]\right)x^k (1-x)^{n-k}\\
&\le \|B_n f-f\| + \sum_{k=1}^{n-1} x^k (1-x)^{n-k}\\
&\le c\,\omega_\varphi^2(f,n^{-1/2}) + \frac{1}{n} \sum_{k=1}^{n-1} p_{n,k}(x)\\
&\le c\,\omega_\varphi^2(f,n^{-1/2}) + \frac{1}{n}.
\end{split}
\end{equation}

We will show that the simultaneous approximation by $\widetilde{B}_n(f)$ satisfies a similar estimate. Before stating that result, let us note that another integer modification of $B_n f$ possesses actually better properties regarding simultaneous approximation. In it, instead of the integer part 
$[\alpha]$ we use the nearest integer $\la\alpha\ra$ to the real $\alpha$. More precisely, if $\alpha\in\R$ is not the arithmetic mean of two consecutive integers, we set $\la\alpha\ra$ to be the integer at which the minimum $\min_{m\in\Z} |\alpha-m|$ is attained. When $\alpha$ is right in the middle between two consecutive integers, we need to impose a tie-breaking rule. Let $m\in\Z$. There are several options:
\begin{itemize}
	\item Round half up: if $\alpha=m+1/2$, then $\la\alpha\ra:=m+1$;
	\item Round half down: if $\alpha=m+1/2$, then $\la\alpha\ra:=m$;
	\item Round half towards zero: if $\alpha=m+1/2$ and $m\ge 0$, then $\la\alpha\ra:=m$; if $\alpha=m+1/2$ and $m<0$, then 
	$\la\alpha\ra:=m+1$;
	\item Round half away from zero: if $\alpha=m+1/2$ and $m\ge 0$, then $\la\alpha\ra:=m+1$; if $\alpha=m+1/2$ and $m<0$, then 
	$\la\alpha\ra:=m$;
	\item Round half to even: if $\alpha=2m\pm 1/2$, then $\la\alpha\ra:=2m$;
	\item Round half to odd: if $\alpha=(2m+1)\pm 1/2$, then $\la\alpha\ra:=2m+1$;
	\item Random half-rounding: if $\alpha=m+1/2$, then $\la\alpha\ra:=m$, or $\la\alpha\ra:=m+1$ with certain probability, which generally depends on $\alpha$.
\end{itemize}
The results we will prove are valid for any tie-breaking rule listed above, including any mixture of them. We will denote that integer modification of the Bernstein polynomial by $\widehat{B}_n(f)$, that is, we set
\[
\widehat{B}_n(f)(x):=\sum_{k=0}^n \la f\left(\frac{k}{n}\right) \binom{n}{k} \ra x^k (1-x)^{n-k}
\]
for $f\in C[0,1]$ and $x\in [0,1]$.

An argument similar to \eqref{eq[B]} yields
\[
\|\widehat{B}_n(f)-f\|\le c\,\omega_\varphi^2(f,n^{-1/2}) + \frac{1}{2n}
\] 
for all $f\in C[0,1]$ with $f(0),f(1)\in\Z$ and all $n\in\N$.

Let us explicitly note that for any fixed $n\ge 2$ the operator $\widetilde{B}_n:C[0,1]\to C[0,1]$ is not bounded in the sense that there does \emph{not} exist a constant $M$ such that
\[
\|\widetilde{B}_n f\|\le M\,\|f\| \quad \forall\,f\in C[0,1].
\] 
That operator is not continuous either. On the other hand, $\widehat{B}_n$ is bounded but not continuous. Both operators are not linear. To emphasize the latter we write $\widetilde{B}_n(f)$ and $\widehat{B}_n(f)$, not $\widetilde{B}_n f$ and $\widehat{B}_n f$.

Recently, we characterized the rate of the simultaneous approximation by the Bernstein operator with Jacobi weights in $L_p$-norm, $1<p\le\infty$, (see \cite{Dr}). In particular, we showed in \cite[Corollary 1.6]{Dr} (with $r=1$) that for all $f\in C^s[0,1]$ and $n\in\N$ there holds
\begin{equation}\label{thmsim}
	\|(B_{n} f)^{(s)} - f^{(s)}\|\\
		\le c\begin{cases}
		\displaystyle{\omega_\varphi^{2}(f',n^{-1/2}) + \omega_1(f',n^{-1})}, &s=1,\\[10pt]
		\displaystyle{\omega_\varphi^{2}(f^{(s)},n^{-1/2}) + \omega_1(f^{(s)},n^{-1})	+ \frac{1}{n}\,\| f^{(s)}\|}, &s\ge 2,
		\end{cases}
\end{equation}
as, moreover, these estimates cannot be improved. Here $\omega_1(F,t)$ is the ordinary modulus of continuity in the uniform norm on the interval $[0,1]$, defined by
\[
\omega_1(F,t):=\sup_{\substack{|x-y|\le t\\x,y\in [0,1]}}|F(x)-F(y)|.
\]

We will verify that the integer forms of the Bernstein polynomials $\widetilde{B}_n$ and $\widehat{B}_n$ satisfy similar direct inequalities. They are stated in the following two theorems.

\begin{thm}\label{thmsim1}
Let $s\in\mathbb{N}_+$. Let $f\in C^s[0,1]$ as $f(0),f(1),f'(0),f'(1)\in\Z$ and $f^{(i)}(0)=f^{(i)}(1)=0$, $i=2,\dots,s$. Let also there exist $n_0\in\N_+$, $n_0\ge s$, such that
\begin{align*}
	f\left(\frac{k}{n}\right) &\ge f(0)+\frac{k}{n}\,f'(0),\quad k=1,\dotsc,s,\ n\ge n_0,\\
	f\left(\frac{k}{n}\right) &\ge f(1)-\left(1-\frac{k}{n}\right)f'(1),\quad k=n-s,\dotsc,n-1,\ n\ge n_0.
\end{align*}
Then for $n\ge n_0$ there holds
\begin{equation*}
\|(\widetilde{B}_n (f))^{(s)}-f^{(s)}\|\le c\begin{cases}
		\displaystyle{\omega_\varphi^{2}(f',n^{-1/2}) + \omega_1(f',n^{-1})}+\frac{1}{n}, &s=1,\\[10pt]
		\displaystyle{\omega_\varphi^{2}(f^{(s)},n^{-1/2}) + \omega_1(f^{(s)},n^{-1})	+ \frac{1}{n}\,\| f^{(s)}\|+\frac{1}{n}}, &s\ge 2.
		\end{cases}
\end{equation*}
The constant $c$ is independent of $f$ and $n$.
\end{thm}

The estimates of the rate of convergence for $\widehat B_n$ are valid under \emph{weaker} assumptions.

\begin{thm}\label{thmsim2}
Let $s\in\mathbb{N}_+$. Let $f\in C^s[0,1]$ as $f(0),f(1),f'(0),f'(1)\in\Z$ and $f^{(i)}(0)=f^{(i)}(1)=0$, $i=2,\dots,s$. Then
\begin{equation*}
\|(\widehat{B}_n (f))^{(s)}-f^{(s)}\|\le c\begin{cases}
		\displaystyle{\omega_\varphi^{2}(f',n^{-1/2}) + \omega_1(f',n^{-1})}+\frac{1}{n}, &s=1,\\[10pt]
		\displaystyle{\omega_\varphi^{2}(f^{(s)},n^{-1/2}) + \omega_1(f^{(s)},n^{-1})	+ \frac{1}{n}\,\|f^{(s)}\|+\frac{1}{n}}, &s\ge 2.
		\end{cases}
\end{equation*}
The constant $c$ is independent of $f$ and $n$.
\end{thm}

We will also show that the assumptions made in Theorems \ref{thmsim1} and \ref{thmsim2} are necessary in order to have uniform simultaneous approximation. The difference between the set of conditions for $s=1$ and $s\ge 2$ is related to the fact that $\widetilde B_n$ and $\widehat B_n$ preserve the polynomials of the form $p\,x+q$, where $p,q\in\Z$. That is verified just as for the Bernstein operators. 

There is an extensive literature on the approximation of functions by polynomials with integer coefficients. A quite helpful introduction to the subject is the monograph \cite{Fe} (see also \cite[Chapter 2, \S\,4]{Lo-Go-Ma:CA}). In particular, the extension of the classical results on simultaneous approximation by algebraic polynomials with real coefficients to the integer case is due to Gelfond \cite{Ge} and Trigub \cite{Tr1,Tr2}. Martinez \cite{FLMa} considered approximation of the derivatives of smooth functions by means of integer forms of the Bernstein polynomials but the coefficients are replaced by their integral part \emph{after} differentiating the Bernstein polynomial of the function.

Finally, let us note that the approximation by polynomials with integer coefficients is important because of their computer implementations.

\section{Proof of the estimates of the rate of convergence}

The integer modifications of the Bernstein polynomials $\widetilde B_n$ and $\widehat B_n$ are not linear. That is why the simplest way to estimate their rate of approximation is to consider their deviation from the linear operator $B_n$ (see \eqref{eq[B]}). We will apply that approach to estimate their rate of simultaneous approximation.

For $n\in\N_+$ and $k=0,\dotsc,n$. We set
\begin{align*}
	\tilde b_n(k)&:=\left[f\left(\frac{k}{n}\right)\binom{n}{k} \right]\,\binom{n}{k}^{-1}
	\intertext{and}
	\hat b_n(k)&:=\la f\left(\frac{k}{n}\right)\binom{n}{k} \ra\,\binom{n}{k}^{-1}.
\end{align*} 
Then the operators $\widetilde{B}_n$ and $\widehat{B}_n$ can be written respectively in the form
\begin{align*}
	\widehat{B}_n (f,x) &=\sum_{k=0}^n \tilde b_n(k)\,p_{n,k}(x)
	\intertext{and}
	\widehat{B}_n (f,x) &=\sum_{k=0}^n \hat b_n(k)\,p_{n,k}(x).
\end{align*}

We will use the forward finite difference operator $\Delta_h$ with step $h$, defined by
\[
\Delta_h f(x):=f(x+h)-f(x),\quad \Delta_h^s:=\Delta_h(\Delta_h^{s-1}).
\]
Then
\begin{equation}\label{eqDelta}
\Delta_h^s f(x)=\sum_{i=0}^s (-1)^i \binom{s}{i} f(x+(s-i)h),\quad x\in [0,1-sh].
\end{equation}
If $h=1$, we will omit the subscript, writing $\Delta:=\Delta_1$. Thus
\begin{equation}\label{eqDelta_b}
\Delta^s \tilde b_n(k)=\sum_{i=0}^s (-1)^i \binom{s}{i} \tilde b_n(k+s-i),\quad k=0,\dotsc,n-s;
\end{equation}
and analogously for $\hat b_n$.

As is known, for $n\ge s$ we have (see \cite{Ma}, or \cite[Chapter 10, (2.3)]{De-Lo:CA}, or \cite[p.\ 125]{Di-To:Mod}) that
\begin{equation}\label{eq1}
(B_n f)^{(s)}(x)=\frac{n!}{(n-s)!} \sum_{k=0}^{n-s} \Delta_{1/n}^{s} f\left(\frac{k}{n}\right) p_{n-s,k}(x), \quad x\in [0,1].
\end{equation}

Similarly, for $n\ge s$ we have
\begin{align}
(\widetilde{B}_n (f))^{(s)}(x)&=\frac{n!}{(n-s)!} \sum_{k=0}^{n-s} \Delta^{s} \tilde b_n(k)\,p_{n-s,k}(x), \quad x\in [0,1],\label{eq2a}
\intertext{and}
(\widehat{B}_n (f))^{(s)}(x)&=\frac{n!}{(n-s)!} \sum_{k=0}^{n-s} \Delta^{s} \hat b_n(k)\,p_{n-s,k}(x), \quad x\in [0,1].\label{eq2}
\end{align}

We proceed to the results that relate $\widetilde{B}_n$ and $\widehat B_n$ to $B_n$.

\begin{thm}\label{tB-Bcl}
Let $s\in\mathbb{N}_+$. Let $f\in C^s[0,1]$ as $f(0),f(1),f'(0),f'(1)\in\Z$ and $f^{(i)}(0)=f^{(i)}(1)=0$, $i=2,\dots,s$. Let also there exist $n_0\in\N_+$, $n_0\ge s$, such that
\begin{align}
	f\left(\frac{k}{n}\right) &\ge f(0)+\frac{k}{n}\,f'(0),\quad k=1,\dotsc,s,\ n\ge n_0,\label{eq27a}\\
	f\left(\frac{k}{n}\right) &\ge f(1)-\left(1-\frac{k}{n}\right)f'(1),\quad k=n-s,\dotsc,n-1,\ n\ge n_0.\label{eq27b}
\end{align}
Then
\[
\|(B_n f)^{(s)}-(\widetilde{B}_n(f))^{(s)}\|\le c\left(\omega_1( f^{(s)},n^{-1}) + \frac{1}{n}\right),\quad n\ge n_0.
\]
The constant $c$ is independent of $f$ and $n$.
\end{thm}

\begin{rem}
Certainly, it suffices to assume instead of the cumbersome \eqref{eq27a}-\eqref{eq27b} that there exists $\delta\in (0,1)$ such that
\begin{align*}
	f(x) &\ge f(0)+x\,f'(0),\quad x\in [0,\delta],\\
	f(x) &\ge f(1)-(1-x)f'(1),\quad x\in [1-\delta,1].
\end{align*}
However, it turns out that the conditions \eqref{eq27a}-\eqref{eq27b} are also necessary unlike the ones above (see \thmref{thmnec2}).
\end{rem}

\begin{thm}\label{tB-B}
Let $s\in\mathbb{N}_+$. Let $f\in C^s[0,1]$ as $f(0),f(1),f'(0),f'(1)\in\Z$ and $f^{(i)}(0)=f^{(i)}(1)=0$, $i=2,\dots,s$. Then
\[
\|(B_n f)^{(s)}-(\widehat{B}_n(f))^{(s)}\|\le c\left(\omega_1( f^{(s)},n^{-1}) + \frac{1}{n}\right).
\]
The constant $c$ is independent of $f$ and $n$.
\end{thm}

Now, Theorems \ref{thmsim1} and \ref{thmsim2} follow directly from \eqref{thmsim} and Theorems \ref{tB-Bcl} and \ref{tB-B}, respectively.

Let us establish Theorems \ref{tB-Bcl} and \ref{tB-B}.

\begin{proof}[Proof of \thmref{tB-Bcl}]
Let $n\ge n_0$. We make use of \eqref{eq1}, \eqref{eq2a}, and the identities $\sum_{j=0}^s \binom{s}{j}=2^s$ and $\sum_{k=0}^{n-s}p_{n-s,k}(x)\equiv 1$ to get
\begin{multline}\label{eq3a}
\left|(B_n f)^{(s)}(x)-(\widetilde{B}_n (f))^{(s)}(x)\right|\\
 \le 2^s\,n^s \max_{0\le k\le n} \left(f\left(\frac{k}{n}\right)-\tilde b_n(k)\right), \quad x\in [0,1].
\end{multline}
Note that $f(k/n)-\tilde b_n(k)\ge 0$, $k=0,\dotsc,n$, because $[\alpha]\le\alpha$.

We will estimate $f(k/n)-\tilde b_n(k)$ separately for $k\le s$, $s+1\le k\le n-s-1$, and $k\ge n-s$. For the middle part, we simply use that if $n\ge 2s+2$, then
\begin{equation}\label{eq4a}
\begin{split}
	f\left(\frac{k}{n}\right)-\tilde b_n(k)
	&= \left(f\left(\frac{k}{n}\right)\binom{n}{k}-\left[f\left(\frac{k}{n}\right)\binom{n}{k} \right]\right)\binom{n}{k}^{-1}\\
	&\le \binom{n}{s+1}^{-1}\le \frac{c}{n^{s+1}},\quad k=s+1,\dotsc,n-s-1.
\end{split}
\end{equation}

Next, we will show that
\begin{equation}\label{eq5a}
f\left(\frac{k}{n}\right)-\tilde b_n(k)
\le \frac{c}{n^s}\,\omega_1( f^{(s)},n^{-1}),\quad k=0,\dotsc,s.
\end{equation}
We apply Taylor's formula, as we take into consideration that $f^{(i)}(0)=0$ for $i=2,\dotsc,s$, to arrive at
\begin{multline}\label{eqT}
f\left(\frac{k}{n}\right)
=f(0) + \frac{k}{n}\,f'(0)\\
+\frac{1}{(s-1)!}\int_0^{k/n} \left(\frac{k}{n}-t\right)^{s-1} \left(f^{(s)}(t)-f^{(s)}(0)\right)dt.
\end{multline}
That implies
\begin{equation}\label{eq7a}
\begin{split}
	f\left(\frac{k}{n}\right)-\left(f(0) + \frac{k}{n}\,f'(0)\right)
&\le \frac{1}{s!}\left(\frac{k}{n}\right)^s \omega_1\left(f^{(s)},\frac{k}{n}\right)\\
&\le \frac{c}{n^s}\,\omega_1( f^{(s)},n^{-1}),\quad k=0,\dotsc,s.
\end{split}
\end{equation}
At the second estimate, we have taken into account the well-known property of the modulus of continuity
\[
\omega_1(F,rt)\le r\omega_1(F,t),
\]
where $r\in\N_+$.

On the other hand, \eqref{eq27a} and
\begin{equation}\label{eq29}
f(0)\binom{n}{k} + f'(0)\,\frac{k}{n}\,\binom{n}{k}\in\Z,
\end{equation}
imply
$$
\left[f\left(\frac{k}{n}\right)\binom{n}{k}\right] \ge f(0)\binom{n}{k} + f'(0)\,\frac{k}{n}\,\binom{n}{k},\quad k=0,\dotsc,s.
$$
Consequently,
\begin{equation}\label{eq9a}
\tilde b_n(k) \ge f(0) + \frac{k}{n}\,f'(0),\quad k=0,\dotsc,s.
\end{equation}
Estimates \eqref{eq7a} and \eqref{eq9a} imply \eqref{eq5a}.

Finally, we observe that, by symmetry, \eqref{eq5a} yields
\begin{equation}\label{eq6a}
f\left(\frac{k}{n}\right)-\tilde b_n(k)
\le \frac{c}{n^s}\,\omega_1( f^{(s)},n^{-1}),
\quad k=n-s,\dotsc,n.
\end{equation}
More precisely, with $\bar f(x):=f(1-x)$ and
\[
\bar{\tilde b}_n(k):=\left[\bar f\left(\frac{k}{n}\right)\binom{n}{k} \right]\,\binom{n}{k}^{-1}
\] 
we have 
\begin{equation}\label{eq21a}
\begin{split}
	&\bar f\left(\frac{k}{n}\right) = f\left(\frac{n-k}{n}\right),\\ 
	&\bar{\tilde b}_n(k) = \tilde b_n(n-k),\\ 
	&\omega_1(\bar f^{(s)},t) = \omega_1(f^{(s)},t).
\end{split}
\end{equation}
Note also that $\bar f\in C^s[0,1]$, $\bar f(0),\bar f'(0)\in\Z$, $\bar f^{(i)}(0)=0$, $i=2,\dots,s$, and for $n\ge n_0$ and $k=1,\dots,s$ we have by \eqref{eq27b}
\[
\bar f\left(\frac{k}{n}\right)=f\left(\frac{n-k}{n}\right)\ge f(1)-\frac{k}{n}f'(1)=\bar f(0)+\frac{k}{n}\bar f'(0).
\]
So, $\bar f$ satisfies the condition \eqref{eq27a} and, in virtue of \eqref{eq5a}, we have
\[
	\bar f\left(\frac{k}{n}\right) -\bar{\tilde b}_n(k)
	\le \frac{c}{n^s}\,\omega_1(\bar f^{(s)},n^{-1}),\quad k=0,\dotsc,s.
\]
As we take into account \eqref{eq21a}, we get \eqref{eq6a}.

Inequalities \eqref{eq3a}-\eqref{eq5a} and \eqref{eq6a} imply the assertion of the theorem.
\end{proof}

We will use the following elementary lemma in the proof the theorem about $\widehat B_n$.

\begin{lem}\label{lm}
Let $m\in\Z$ and $\alpha,\omega\in\R$. If $|\alpha-m|\le\omega$, then $|\la\alpha\ra-m|\le 2\omega$.
\end{lem}

\begin{proof}
If $\omega<1/2$, then $\la\alpha\ra=m$. If, on the other hand, $\omega\ge 1/2$, then
\[
|\la\alpha\ra-m|\le |\la\alpha\ra-\alpha|+|\alpha-m|\le \frac{1}{2}+\omega\le 2\omega.
\]
\end{proof}

\begin{proof}[Proof of \thmref{tB-B}]
We proceed similarly to the proof of the previous theorem. Since the assertion is trivial for $n<s$, we assume that $n\ge s$. We make use of \eqref{eq1} and \eqref{eq2} to get
\begin{multline}\label{eq3}
\left|(B_n f)^{(s)}(x)-(\widehat{B}_n (f))^{(s)}(x)\right|\\
\le 2^s\,n^s \max_{0\le k\le n} \left|f\left(\frac{k}{n}\right)-\hat b_n(k)\right|, \quad x\in [0,1].
\end{multline}

Again we estimate separately the terms $|f(k/n)-\hat b_n(k)|$ for $k\le s$, $s+1\le k\le n-s-1$, and $k\ge n-s$. For the middle part, we have similarly to \eqref{eq4a}
\begin{equation}\label{eq4}
\left|f\left(\frac{k}{n}\right)-\hat b_n(k)\right|\le \frac{c}{n^{s+1}},\quad k=s+1,\dotsc,n-s-1,\ n\ge 2s+2.
\end{equation}

Next, we will show that
\begin{equation}\label{eq5}
\left|f\left(\frac{k}{n}\right) -\hat b_n(k)\right|\le \frac{c}{n^s}\,\omega_1( f^{(s)},n^{-1}),\quad k=0,\dotsc,s.
\end{equation}
In virtue of \eqref{eqT}, we have
\begin{equation}\label{eq7}
\left|f\left(\frac{k}{n}\right) - \left(f(0) + \frac{k}{n}\,f'(0)\right)\right|
\le \frac{c}{n^s}\,\omega_1( f^{(s)},n^{-1}),\quad k=0,\dotsc,s.
\end{equation}
That implies
\begin{multline}\label{eq8}
\left|f\left(\frac{k}{n}\right)\binom{n}{k} - \left(f(0)\binom{n}{k} + f'(0)\,\frac{k}{n}\,\binom{n}{k}\right)\right|\\
\le \frac{c}{n^s}\,\binom{n}{k}\omega_1( f^{(s)},n^{-1}),\quad k=0,\dotsc,s.
\end{multline}
We apply \lemref{lm} with
\begin{align*}
	\alpha &=f\left(\frac{k}{n}\right)\binom{n}{k},\\
	m &=f(0)\binom{n}{k} + f'(0)\,\frac{k}{n}\,\binom{n}{k}\in\Z,\\
	\omega &=\frac{c}{n^s}\,\binom{n}{k}\omega_1( f^{(s)},n^{-1}),
\end{align*}
where the constant $c$ is the one on the right-hand side of \eqref{eq8}.

Thus we arrive at
\begin{multline*}
\left|\la f\left(\frac{k}{n}\right)\binom{n}{k}\ra - \left(f(0)\binom{n}{k} + f'(0)\,\frac{k}{n}\,\binom{n}{k}\right)\right|\\
\le \frac{c}{n^s}\,\binom{n}{k}\omega_1( f^{(s)},n^{-1}),\quad k=0,\dotsc,s,
\end{multline*}
and, consequently,
\begin{equation}\label{eq9}
\left|\hat b_n(k) - \left(f(0) + \frac{k}{n}\,f'(0)\right)\right|
\le \frac{c}{n^s}\,\omega_1( f^{(s)},n^{-1}),\quad k=0,\dotsc,s.
\end{equation}
Estimates \eqref{eq7} and \eqref{eq9} yield \eqref{eq5}.

Finally, we derive
\begin{equation}\label{eq6}
\left|f\left(\frac{k}{n}\right) - \hat b_n(k)\right|\le \frac{c}{n^s}\,\omega_1( f^{(s)},n^{-1}),
\quad k=n-s,\dotsc,n.
\end{equation}
from \eqref{eq5} by symmetry just as in the proof of \eqref{eq6a} with $\bar{\tilde b}_n(k)$ replaced with
\[
\bar{\hat b}_n(k):=\la \bar f\left(\frac{k}{n}\right)\binom{n}{k} \ra\,\binom{n}{k}^{-1}.
\] 

Inequalities \eqref{eq3}-\eqref{eq5} and \eqref{eq6} imply the assertion of the theorem.
\end{proof}

\section{Optimality of the assumptions in Theorems \ref{thmsim1} and \ref{thmsim2}}

We will establish the necessity of the assumptions made in Theorems \ref{thmsim1} and \ref{thmsim2}. We begin with the operator $\widehat B_n$ since stronger results are valid for it.

First of all, let us note that if 
\begin{equation}\label{eq23a}
\lim_{n\to\infty}\|\widehat{B}_n(f)-f\|=0\quad\text{and}\quad
\lim_{n\to\infty}\|(\widehat{B}_n(f))^{(s)}-f^{(s)}\|=0,
\end{equation}
then $f^{(i)}(0),f^{(i)}(1)\in\Z$ for $i=0,\dotsc,s$. Indeed, as is known for any $g\in C^s[0,1]$ we have (see e.g.~\cite[Chapter 2, Theorem 5.6]{De-Lo:CA})
\[
\|g^{(i)}\|\le c\big(\|g\| + \|g^{(s)}\|\big),\quad i=1,\dotsc,s-1.
\]
Therefore \eqref{eq23a} implies
\begin{equation}\label{eq11}
\lim_{n\to\infty}\|(\widehat{B}_n(f))^{(i)}-f^{(i)}\|=0,\quad i=0,\dotsc,s;
\end{equation}
hence $f^{(i)}(0),f^{(i)}(1)\in\Z$ for $i=0,\dotsc,s$. A similar result holds for $\widetilde{B}_n$.

\begin{thm}\label{thmnec1}
Let $s\in\N_+$, as $s\ge 2$, and $f\in C^s[0,1]$. If
\begin{equation}\label{eq23}
\lim_{n\to\infty}\|\widehat{B}_n(f)-f\|=0\quad\text{and}\quad
\lim_{n\to\infty}\|(\widehat{B}_n(f))^{(s)}-f^{(s)}\|=0,
\end{equation}
then 
$f^{(i)}(0)=f^{(i)}(1)=0$, $i=2,\dotsc,s$.
\end{thm}

\begin{proof}
It is sufficient to establish the theorem at the point $x=0$; for $x=1$ it follows by symmetry. We use induction on $s$.

Let $s=2$. Relation \eqref{eq11}, in particular, yields
\[
\lim_{n\to\infty}(\widehat B_n(f))'(0)=f'(0),
\] 
that is (see \eqref{eq2} with $s=1$),
\begin{equation}\label{eq12}
\lim_{n\to\infty} n\Delta \hat b_n(0)=f'(0).
\end{equation}
Since $n\Delta \hat b_n(0)\in\Z$ for all $n$, \eqref{eq12} implies
\[
n\Delta \hat b_n(0)=f'(0)\quad\text{for $n$ large enough};
\]
hence
\begin{equation}\label{eq13}
\hat b_n(1)=\hat b_n(0)+\frac{1}{n}\,f'(0)=f(0)+\frac{1}{n}\,f'(0)
\end{equation}
Similarly, from $\lim_{n\to\infty}(\widehat B_n(f))''(0)=f''(0)$ we derive
\begin{equation}\label{eq14}
n(n-1)\Delta^2 \hat b_n(0)=f''(0)\quad\text{for $n$ large enough}.
\end{equation}
By Taylor's formula, we have
\begin{equation}\label{eq38}
f\left(\frac{2}{n}\right)=f(0)+\frac{2}{n}\,f'(0) + \frac{2}{n^2}\,f''(0) + \int_0^{2/n} \left(\frac{2}{n}-t\right)(f''(t)-f''(0))\,dt.
\end{equation}
Next, we proceed similarly to the proof of \thmref{tB-B}. We multiply both sides of the above identity by $\binom{n}{2}$ and rearrange the terms to get
\begin{multline}\label{eq27}
	f\left(\frac{2}{n}\right)\binom{n}{2}-\left(f(0)\binom{n}{2} + (n-1)f'(0) + f''(0)\right)\\
	=-\frac{1}{n}\,f''(0)+\binom{n}{2}\int_0^{2/n} \left(\frac{2}{n}-t\right)(f''(t)-f''(0))\,dt.
\end{multline}
Consequently,
\begin{multline*}
	\left|f\left(\frac{2}{n}\right)\binom{n}{2}-\left(f(0)\binom{n}{2} + (n-1)f'(0) + f''(0)\right)\right|\\
	\le \frac{1}{n}\,|f''(0)|+\omega_1\left(f'',\frac{2}{n}\right),
\end{multline*}
which shows that for large $n$ we have
\[
\la f\left(\frac{2}{n}\right)\binom{n}{2}\ra = f(0)\binom{n}{2} + (n-1)f'(0) + f''(0).
\] 
Therefore
\begin{equation}\label{eq15}
\hat b_n(2)=f(0)+\frac{2}{n}\,f'(0)+\frac{2}{n(n-1)}\,f''(0)\quad\text{for $n$ large enough}.
\end{equation}
Now, fixing some $n$ large enough, we deduce from \eqref{eq13}-\eqref{eq15} that
\begin{align*}
	&f''(0) = n(n-1)(\hat b_n(2)-2\hat b_n(1)+\hat b_n(0))\\
	&\quad =n(n-1)\left(f(0)+\frac{2}{n}\,f'(0)+\frac{2}{n(n-1)}\,f''(0) - 2\left(f(0)+\frac{1}{n}\,f'(0)\right) + f(0)\right)\\
	&\quad =2f''(0);
\end{align*}
hence $f''(0)=0$.

Let the assertion of the theorem hold for some $s-1$, $s\ge 3$. We will prove that then it holds for $s$ too.

As we noted in the beginning of the section, \eqref{eq23} implies
\[
\lim_{n\to\infty}\|(\widehat{B}_n(f))^{(s-1)}-f^{(s-1)}\|=0.
\]
Hence, in virtue of the induction hypothesis, we have $f^{(i)}(0)=0$ for $i=2,\dotsc,s-1$.

By Taylor's formula we have
\begin{multline}\label{eq16}
	f\left(\frac{k}{n}\right)
	=f(0)+\frac{k}{n}\,f'(0) + \left(\frac{k}{n}\right)^s\frac{f^{(s)}(0)}{s!}\\ 
	+ \frac{1}{(s-1)!}\int_0^{k/n} \left(\frac{k}{n}-t\right)^{s-1} (f^{(s)}(t)-f^{(s)}(0))\,dt.
\end{multline}
We multiply both sides by $\binom{n}{k}$. For $1\le k<s$ we derive the inequality
\begin{align*}
	\bigg|f\left(\frac{k}{n}\right)\binom{n}{k} - &\left(f(0)\binom{n}{k}+f'(0)\frac{k}{n}\binom{n}{k}\right)\bigg|\\
	&\qquad \le \binom{n}{k}\left(\frac{k}{n}\right)^s\frac{1}{s!} \left(|f^{(s)}(0)| + \omega_1\left(f^{(s)},\frac{k}{n}\right)\right)\\
	&\qquad \le \frac{c}{n}\left(|f^{(s)}(0)| + \omega_1( f^{(s)},n^{-1})\right).
\end{align*}
Consequently, for large $n$ we have 
\[
\la f\left(\frac{k}{n}\right)\binom{n}{k}\ra = f(0)\binom{n}{k}+f'(0)\frac{k}{n}\binom{n}{k};
\]
hence
\begin{equation}\label{eq17}
\hat b_n(k)=f(0)+\frac{k}{n}\,f'(0)\quad\text{for $0\le k<s$ and large $n$}.
\end{equation}

In order to calculate $\hat b_n(s)$, we observe that
\begin{equation*}
\lim_{n\to\infty} \binom{n}{s}\left(\frac{s}{n}\right)^s = \frac{s^s}{s!}.
\end{equation*}
We proceed just as in this case $s=2$: we multiple both sides of \eqref{eq16} by $\binom{n}{s}$ and rearrange the terms to arrive at 
\begin{align*}
	\bigg|f\left(\frac{s}{n}\right)\binom{n}{s} &- \left(f(0)\binom{n}{s}+f'(0)\frac{s}{n}\binom{n}{s} + \frac{s^s}{(s!)^2}\,f^{(s)}(0)\right)\bigg|\\
	&\qquad \le \left(\frac{s^s}{s!}-\binom{n}{s}\left(\frac{s}{n}\right)^s\right)\frac{1}{s!}\,|f^{(s)}(0)| 
	+ \frac{1}{s!}\binom{n}{s}\left(\frac{s}{n}\right)^s\omega_1\left(f^{(s)},\frac{s}{n}\right)\\
	&\qquad \le \frac{c}{n}|f^{(s)}(0)| + c\,\omega_1( f^{(s)},n^{-1}).
\end{align*}
Consequently, for large $n$
\[
\la f\left(\frac{s}{n}\right)\binom{n}{s} \ra = f(0)\binom{n}{s}+f'(0)\frac{s}{n}\binom{n}{s} + \la \frac{s^s}{(s!)^2}\,f^{(s)}(0)\ra + r_s,
\]
where $r_s\in\{-1,0,1\}$. Consequently,
\begin{equation}\label{eq18}
\hat b_n(s)=f(0)+\frac{s}{n}\,f'(0)+\left(\la \frac{s^s}{(s!)^2}\,f^{(s)}(0)\ra + r_s\right)\binom{n}{s}^{-1}.
\end{equation}

Relations \eqref{eq17} and \eqref{eq18} yield
\begin{equation}\label{eq19}
\Delta^s \hat b_n(0)=\left(\la \frac{s^s}{(s!)^2}\,f^{(s)}(0)\ra + r_s\right)\binom{n}{s}^{-1}.
\end{equation}

On the other hand, since $\lim_{n\to\infty}\|(\widehat{B}_n(f))^{(s)}-f^{(s)}\|=0$, and, in particular, 
$\lim_{n\to\infty}(\widehat{B}_n(f))^{(s)}(0)=f^{(s)}(0)$, we have that
\[
\lim_{n\to\infty} \frac{n!}{(n-s)!} \Delta^s \hat b_n(0)=f^{(s)}(0).
\]
Taking into account that
\[
\frac{n!}{(n-s)!} \Delta^s \hat b_n(0)\in\Z\quad\forall n,
\]
we deduce that for large $n$ there holds
\[
\frac{n!}{(n-s)!} \Delta^s \hat b_n(0)=f^{(s)}(0).
\]
That, in combination with \eqref{eq19}, yields
\begin{equation}\label{eq20}
s! \left(\la \frac{s^s}{(s!)^2}\,f^{(s)}(0)\ra + r_s\right)=f^{(s)}(0)\quad\text{for $n$ large enough}.
\end{equation}
First of all, this relation implies that the integer $f^{(s)}(0)$ is divisible by $s!$, i.e.~$f^{(s)}(0)=s!\,m_s$ with some $m_s\in\Z$. Thus
\eqref{eq20} can be reduced to
\[
\la \frac{s^s}{s!}\,m_s\ra + r_s = m_s.
\]
Consequently,
\[
|m_s|\left(\frac{s^s}{s!}-1\right)\le \frac{3}{2}.
\]
It remains to take into account that $s^s/s!$ increases on $s$; hence $s^s/s!\ge 9/2$ for $s\ge 3$, and then $|m_s|\le 3/7$, which is possible only if $m_s=0$. Thus $f^{(s)}(0)=0$.
\end{proof}

Necessary conditions for the simultaneous approximation by means of $\widetilde B_n$ are given in the following theorem.

\begin{thm}\label{thmnec2}
Let $s\in\N_+$ and $f\in C^s[0,1]$. If 
\begin{equation}\label{eq22}
	\lim_{n\to\infty}\|\widetilde{B}_n(f)-f\|=0\quad\text{and}\quad
	\lim_{n\to\infty}\|(\widetilde{B}_n(f))^{(s)}-f^{(s)}\|=0,
\end{equation}
then $f^{(i)}(0)=f^{(i)}(1)=0$, $i=2,\dotsc,s$, and there exists $n_0\in\N_+$, $n_0\ge s$, such that
\begin{align}
	f\left(\frac{k}{n}\right) &\ge f(0)+\frac{k}{n}\,f'(0),\quad k=1,\dotsc,s,\ n\ge n_0,\label{eq22a}\\
	f\left(\frac{k}{n}\right) &\ge f(1)-\left(1-\frac{k}{n}\right)f'(1),\quad k=n-s,\dotsc,n-1,\ n\ge n_0.\notag
\end{align}
\end{thm}

\begin{proof}
It is sufficient to establish the theorem at the point $x=0$; for $x=1$ it follows by symmetry.

We argue as in the proof of the preceding theorem. However, here more efforts are required.

Using induction on $s$, we will prove that $f^{(i)}(0)=0$, $i=2,\dotsc,s$ and
\begin{equation}\label{eq25}
\tilde b_n(k)=f(0)+\frac{k}{n}\,f'(0), \quad k=1,\dotsc,s,\ n\ge n_0.
\end{equation}
with some $n_0$. The latter implies directly the inequalities \eqref{eq22a} because
\[
f\left(\frac{k}{n}\right)\ge \left[f\left(\frac{k}{n}\right)\binom{n}{k}\right]\binom{n}{k}^{-1}
=f(0)+\frac{k}{n}\,f'(0), \quad k=1,\dotsc,s,\ n\ge n_0.
\]

As in the proof of \thmref{thmnec1}, we deduce from
\[
\lim_{n\to\infty}\|(\widetilde{B}_n(f))^{(s)}-f^{(s)}\|=0
\]
that there exists $n_0\in\N$, $n_0\ge s$, such that
\begin{equation}\label{eq28}
\frac{n!}{(n-i)!}\Delta^i \tilde b_n(0)=f^{(i)}(0), \quad i=1,\dotsc,s,\ n\ge n_0.
\end{equation}
That directly yields \eqref{eq25} for $s=1$ and the assertion of the theorem is verified for $s=1$.

In order to complete the proof for larger $s$, we use that if $f\in C^s[0,1]$ and $\lim_{n\to\infty}\|(\widetilde{B}_n(f))^{(s)}-f^{(s)}\|=0$, then
\begin{equation*}
\lim_{n\to\infty}\|(B_n f)^{(s)}-(\widetilde{B}_n(f))^{(s)}\|=0;
\end{equation*}
hence
\begin{equation}\label{eq33}
\lim_{n\to\infty}\left((B_n f)^{(s)}\left(\frac{y}{n}\right)-(\widetilde{B}_n(f))^{(s)}\left(\frac{y}{n}\right)\right)=0,\quad y\in [0,1].
\end{equation}

By \eqref{eqDelta}-\eqref{eq2a}, after reordering the terms, we arrive at the identity
\begin{multline}\label{eq31}
	(B_n f)^{(s)}(x)-(\widetilde{B}_n(f))^{(s)}(x)\\
	=\frac{n!}{(n-s)!} 
	\sum_{k=0}^{n-s}\sum_{j=k}^{k+s} (-1)^{s+j-k}\binom{s}{j-k}\left(f\left(\frac{j}{n}\right)-\tilde b_n(j)\right)p_{n-s,k}(x).
\end{multline}

We observe that, in virtue of \eqref{eq4a}, for $n\ge 3s+2$ and $x\in [0,1]$ there holds (cf.\ \eqref{eq3a})
\begin{align}
	&\left|\sum_{k=s+1}^{n-2s-1}\sum_{j=k}^{k+s} 
	(-1)^{s+j-k}\binom{s}{j-k}\left(f\left(\frac{j}{n}\right)-\tilde b_n(j)\right)p_{n-s,k}(x)\right|
	\le\frac{c}{n^{s+1}}\label{eq34}
	\intertext{and}
	&\left|\sum_{k=1}^{s}\sum_{j=s+1}^{k+s} 
	(-1)^{s+j-k}\binom{s}{j-k}\left(f\left(\frac{j}{n}\right)-\tilde b_n(j)\right)p_{n-s,k}(x)\right|
	\le\frac{c}{n^{s+1}}.\label{eq34a}
\end{align}

Next, we observe that if $n\ge 4s+1$, then $p_{n-s,k}(y/n)\le c\,n^{-s-1}$ for all $y\in [0,1]$ and $k=n-2s,\dotsc,n-s$. Therefore, taking also into account that $0\le f(j/n)-\tilde b_n(j)\le 1$ (see \eqref{eq4a}), we arrive at
\begin{multline}\label{eq35}
	\left|\sum_{k=n-2s}^{n-s}\sum_{j=k}^{k+s} 
	(-1)^{s+j-k}\binom{s}{j-k}\left(f\left(\frac{j}{n}\right)-\tilde b_n(j)\right)p_{n-s,k}\left(\frac{y}{n}\right)\right|\\
	\le\frac{c}{n^{s+1}},\quad y\in [0,1].
\end{multline}

We apply \eqref{eq33}-\eqref{eq35}, reorder the terms and take into account that $\tilde b_n(0)=f(0)$. Thus, for $y\in [0,1]$, we deduce that
\begin{equation}\label{eq36}
	\lim_{n\to\infty} \frac{n!}{(n-s)!}\sum_{j=1}^{s}
	(-1)^{s-j}\left(f\left(\frac{j}{n}\right)-\tilde b_n(j)\right)\sum_{k=0}^j (-1)^k\binom{s}{j-k}p_{n-s,k}\left(\frac{y}{n}\right)
	=0.
\end{equation}

We will evaluate that limit in another way. Clearly,
\begin{equation}\label{eq37}
\lim_{n\to\infty} \sum_{k=0}^j (-1)^k\binom{s}{j-k} p_{n-s,k}\left(\frac{y}{n}\right)
=\frac{1}{e^y}\sum_{k=0}^j (-1)^k \frac{y^k}{k!}\binom{s}{j-k}.
\end{equation}

Let $s=2$. In virtue of Taylor's formula, we have
\[
f\left(\frac{j}{n}\right)=f(0)+\frac{j}{n}\,f'(0) + \frac{j^2}{n^2}\,\frac{f''(0)}{2} + o(n^{-2}),\quad j=1,2.
\]
Taking into account \eqref{eq25} with $k=s=1$ and \eqref{eq28} with $i=s=2$, we derive
\begin{equation}\label{eqb2}
\tilde b_n(2)=f(0)+\frac{2}{n}\,f'(0)+\frac{1}{n(n-1)}\,f''(0),\quad n\ge n_0.
\end{equation}
Consequently,
\begin{align}
	&n(n-1)\left(f\left(\frac{1}{n}\right)-\tilde b_n(1)\right)=\frac{f''(0)}{2}\, + o(1)
	\intertext{and}
	&n(n-1)\left(f\left(\frac{2}{n}\right)-\tilde b_n(2)\right)=f''(0) + o(1)
\end{align}
Now, we substitute the last two relations into \eqref{eq36} with $s=2$ and $y=1$ and take into account \eqref{eq37} with $y=1$, to deduce that $f''(0)=0$; hence, in virtue of \eqref{eqb2}, we also get \eqref{eq25} for $s=2$. Thus, the assertion of the theorem is verified for $s=2$. 

We proceed by induction on $s$. Relations \eqref{eq22} imply 
\[
\lim_{n\to\infty}\|(\widetilde{B}_n(f))^{(s-1)}-f^{(s-1)}\|=0.
\]
Therefore, in virtue of the induction hypothesis, we have that $f^{(i)}(0)=0$, $i=2,\dotsc,s-1$, $s\ge 3$, and
\begin{equation}\label{eq25a}
\tilde b_n(j)=f(0)+\frac{j}{n}\,f'(0), \quad j=1,\dotsc,s-1,\ n\ge n_0.
\end{equation}
Then Taylor's formula yields
\[
f\left(\frac{j}{n}\right)=f(0)+\frac{j}{n}\,f'(0) + \frac{j^s}{n^s}\,\frac{f^{(s)}(0)}{s!} + o(n^{-s}),\quad j=1,\dotsc,s.
\]
The relations \eqref{eq28} with $i=s$ and \eqref{eq25a} imply
\begin{equation}\label{eq41}
\tilde b_n(s)=f(0)+\frac{s}{n}\,f'(0) + \frac{(n-s)!}{n!}\,f^{(s)}(0),\quad n\ge n_0.
\end{equation}
Therefore
\begin{align*}
	&\frac{n!}{(n-s)!}\left(f\left(\frac{j}{n}\right)-\tilde b_n(j)\right)=\frac{j^s}{s!}\,f^{(s)}(0)+o(1),\quad j=1,\dotsc,s-1,
	\intertext{and}
	&\frac{n!}{(n-s)!}\left(f\left(\frac{s}{n}\right)-\tilde b_n(s)\right)= \left(\frac{s^s}{s!}-1\right)f^{(s)}(0)+ o(1).
\end{align*}

Now, if we substitute the last two relations into \eqref{eq36}, we arrive at
\begin{equation*}
f^{(s)}(0)\sum_{k=0}^{s} \frac{(-1)^k y^k}{k!}\left(\binom{s}{k}-\sum_{j=k}^s (-1)^{s-j}\frac{j^s}{s!}\binom{s}{j-k}\right)=0,\quad y\in [0,1]
\end{equation*}
(actually the summand for $k=0$ is $0$). Consequently, the coefficient of $y^s$ is equal to zero, that is, 
\[
\frac{(-1)^s f^{(s)}(0)}{s!}\left(1-\frac{s^s}{s!}\right)=0.
\]
Therefore $f^{(s)}(0)=0$ and then, in virtue of \eqref{eq41}, $\tilde b_n(s)=f(0)+\frac{s}{n}\,f'(0)$.
\end{proof}

\bigskip
\begin{footnotesize}
\noindent
\begin{tabular}{ll}
Borislav R. Draganov& \\
Dept. of Mathematics and Informatics&
Inst. of Mathematics and Informatics\\
University of Sofia&
Bulgarian Academy of Sciences\\
5 James Bourchier Blvd.&
bl. 8 Acad. G. Bonchev Str.\\
1164 Sofia&
1113 Sofia\\
Bulgaria&
Bulgaria\\
bdraganov@fmi.uni-sofia.bg& \\
\end{tabular}

\end{footnotesize}


\begin{thebibliography}{50}

\bibitem{Fe}
Le Baron O. Ferguson, Approximation by Polynomials with Integral Coefficients, Mathematical Surveys Vol. 17, American Mathematical Society, 1980.

\bibitem{Bu:BP}
J. Bustamante, Bernstein Operators and Their Properties, Birkh\"auser, 2017.

\bibitem{De-Lo:CA}
R. A. DeVore, G. G. Lorentz, Constructive Approximation, Springer-Verlag, Berlin, 1993.

\bibitem{Di-To:Mod}
Z. Ditzian, V. Totik, Moduli of Smoothness, Springer-Verlag, New York, 1987.

\bibitem{Dr-Iv:Char}
B. R. Draganov, K. G. Ivanov, A New characterization of weighted Peetre $K$-functionals, Constr. Approx. 21 (2005), 113--148.

\bibitem{Dr}
B. R. Draganov, Strong estimates of the weighted simultaneous approximation by the Bernstein and Kantorovich operators and their iterated Boolean sums, J. Approx. Theory 200 (2015), 92--135.

\bibitem{Dz-Sh}
V. K. Dzyadyk, I. A. Shevchuk, Theory of Uniform Approximation of Functions by Polynomials, Walter de Gruyter, Berlin, 2008.

\bibitem{Ge}
A. O. Gelfond, On uniform approximation by polynomials with rational integral coefficients, Uspekhi Mat. Nauk 10 (1955), 41--65 (in Russian).

\bibitem{Iv:Dir}
K. G. Ivanov, Some characterizations of the best algebraic approximation in $L_p[-1,1]$ $(1\le p \le \infty)$, C. R. Acad. Bulgare Sci.
34, 1229--1232 (1981).

\bibitem{Iv:Char}
K. G. Ivanov, A characterization of weighted Peetre $K$-functionals, J. Approx. Theory 56 (1989), 185--211.

\bibitem{Ka}
L. V. Kantorovich, Some remarks on the approximation of functions by means of polynomials with integer coefficients, Izv. Akad. Nauk SSSR, Ser. Mat. 9 (1931), 1163--1168 (in Russian).

\bibitem{Ko-Le-Sh:2}
K. Kopotun, D. Leviatan, I. A. Shevchuk, New moduli of smoothness, Publ. Math. Inst. (Beograd) 96 (110) (2014), 169--180.

\bibitem{Ko-Le-Sh}
K. Kopotun, D. Leviatan, I. A. Shevchuk, New moduli of smoothness: weighted DT moduli revisited and applied, Constr. Approx. 42 (2015), 129--159.

\bibitem{Ko-Le-Sh:3}
K. Kopotun, D. Leviatan, I.A. Shevchuk, On weighted approximation with Jacobi weights, 2017, arXiv:1710.05059.

\bibitem{Lo-Go-Ma:CA}
G. G. Lorentz, M. v.Golitschek, Y. Makovoz, Constructive Approximation, Advanced Problems, Springer-Verlag, Berlin, 1996. 

\bibitem{FLMa}
F. L. Martinez, Some properties of two-dimensional Bernstein polynomials, J. Approx. Theory 59 (1989), 300--306.

\bibitem{Ma}
R. Martini, On the approximation of functions together with their derivatives by certain linear positive operators, Indag. Math. 31 (1969), 473--481. 

\bibitem{Tr1}
R. M. Trigub, Approximation of functions by polynomials with integer coefficients, Dokl. Akad. Nauk SSSR 140 (1961), 773--775 (in Russian).

\bibitem{Tr2}
R. M. Trigub, Approximation of functions by polynomials with integer coefficients, Izv. Akad. Nauk SSSR Ser. Mat. 26 (1962), 261--280 (in Russian).

\end{thebibliography}
\end{document}